\documentclass[twocolumn,british]{IEEEtran}
\usepackage[T1]{fontenc}
\usepackage[latin9]{inputenc}
\usepackage{refstyle}
\usepackage{float}
\usepackage{mathtools}
\usepackage{url}
\usepackage{enumitem}
\usepackage{amsmath}
\usepackage{amsthm}
\usepackage{amssymb}
\usepackage{graphicx}

\makeatletter

\AtBeginDocument{\providecommand\secref[1]{\ref{sec:#1}}}
\AtBeginDocument{\providecommand\tabref[1]{\ref{tab:#1}}}
\AtBeginDocument{\providecommand\eqref[1]{\ref{eq:#1}}}
\AtBeginDocument{\providecommand\thmref[1]{\ref{thm:#1}}}
\AtBeginDocument{\providecommand\enuref[1]{\ref{enu:#1}}}
\AtBeginDocument{\providecommand\lemref[1]{\ref{lem:#1}}}
\AtBeginDocument{\providecommand\defref[1]{\ref{def:#1}}}
\AtBeginDocument{\providecommand\propref[1]{\ref{prop:#1}}}
\AtBeginDocument{\providecommand\figref[1]{\ref{fig:#1}}}
\AtBeginDocument{\providecommand\assmref[1]{\ref{assm:#1}}}
\AtBeginDocument{\providecommand\subsecref[1]{\ref{subsec:#1}}}
\providecommand{\tabularnewline}{\\}

\floatstyle{ruled}
\newfloat{algorithm}{tbp}{loa}
\providecommand{\algorithmname}{Algorithm}
\floatname{algorithm}{\protect\algorithmname}
\RS@ifundefined{subsecref}
  {\newref{subsec}{name = \RSsectxt}}
  {}
\RS@ifundefined{thmref}
  {\def\RSthmtxt{theorem~}\newref{thm}{name = \RSthmtxt}}
  {}
\RS@ifundefined{lemref}
  {\def\RSlemtxt{lemma~}\newref{lem}{name = \RSlemtxt}}
  {}

\usepackage[natbibapa]{apacite}
\theoremstyle{plain}
\newtheorem{thm}{\protect\theoremname}
\theoremstyle{definition}
\newtheorem{defn}[thm]{\protect\definitionname}
\theoremstyle{plain}
\newtheorem{lem}[thm]{\protect\lemmaname}
\theoremstyle{plain}
\newtheorem{prop}[thm]{\protect\propositionname}

\usepackage[hidelinks]{hyperref}
\usepackage{amsmath,amssymb,amsfonts}
\usepackage{algpseudocode}
\usepackage{graphicx}
\usepackage{textcomp}
\usepackage{placeins}
\usepackage{xcolor}
\usepackage{setspace}
\usepackage{caption}
\usepackage{microtype}
\usepackage[all]{nowidow}
\usepackage{relsize}
\usepackage{scalefnt}
\usepackage{lmodern}

\newtheorem{assm}{Assumption}
\newref{eq}{refcmd={(\ref{#1})}}
\newref{fig}{refcmd={Fig.~\ref{#1}}}
\newref{def}{refcmd={Definition~\ref{#1}}}
\newref{lem}{refcmd={Lemma~\ref{#1}}}
\newref{sec}{refcmd={Section~\ref{#1}}}
\newref{subsec}{refcmd={Subsection~\ref{#1}}}
\newref{tab}{refcmd={Table~\ref{#1}}}
\newref{thm}{refcmd={Theorem~\ref{#1}}}
\newref{assm}{refcmd={Assumption~\ref{#1}}}
\DeclareMathOperator*{\argmin}{argmin\ }

\providecommand{\keywords}[1]{\textbf{\textit{Index terms---}} #1}

\makeatother

\usepackage{babel}
\providecommand{\definitionname}{Definition}
\providecommand{\lemmaname}{Lemma}
\providecommand{\propositionname}{Proposition}
\providecommand{\theoremname}{Theorem}

\begin{document}

\title{Undiscounted Control Policy Generation for Continuous-Valued Optimal Control by Approximate Dynamic Programming}

\author{Jonathan Lock, Tomas McKelvey}
\maketitle
\begin{abstract}
We present a numerical method for generating the state-feedback control policy associated with general undiscounted, constant-setpoint, infinite-horizon, nonlinear optimal control problems with continuous state variables. The method is based on approximate dynamic programming, and is closely related to approximate policy iteration. Existing methods typically terminate based on the convergence of the control policy and either require a discounted problem formulation or demand the cost function to lie in a specific subclass of functions. The presented method extends on existing termination criterea by requiring both the control policy and the resulting system state to converge, allowing for use with undiscounted cost functions that are bounded and continuous. This paper defines the numerical method, derives the relevant underlying mathematical properties, and validates the numerical method with representative examples. A MATLAB implementation with the shown examples is freely available.
\end{abstract}

\\

\medskip{}

\section{Introduction}

Practical methods for generating the optimal control policy (i.e.~the state feedback function) for general non-linear optimal control problems are useful tools for control engineers. If the optimal control policy is known, a real-time optimal controller can be implemented on very computationally limited hardware as the optimal control signal can be generated simply by interpolating the pre-computed optimal control based on the current system state. However, one practical difficultly lies in pre-computing the optimal control policy, which can be very computationally expensive. Although several methods for solving this class of problem are well-studied, dynamic programming (DP) variants being one example, they all have associated limitations or drawbacks. Policy iteration is one extensively studied variant of DP (e.g.~\citealt[p. 246]{Bertsekas2017}; \citealt[p. 295]{Puterman1994}; \citealt{puterman_convergence_1979}) that has been used for over 40 years for finding the optimal control policy for discrete-valued, non-linear, infinite-horizon problems, i.e.~where the state and control variables are taken from discrete sets. 

Approximate dynamic programming (ADP) is another well-known extension of DP (see for instance \citet{powell_what_2009} for a general introduction) that approximates the cost function using a prescribed set of basis functions. One group of ADP methods approximate the cost function by interpolating costs and optimal controls between discrete gridded points (e.g.~\citet{munos_variable_2002,santos_analysis_1998}). This approach allows for extending DP to applications with continuous state variables.

Assuming the problem of finding the approximately-optimal control policy for continuous-valued, non-linear, infinite-horizon problems, one might attempt to use traditional policy iteration in concert with ADP. However, this is problematic as traditional policy iteration requires the set of states and controls to be discrete (i.e.~finite) to terminate, while the interpolation performed with ADP leads to a continuous (i.e.~infinite) number of possible states and controls. This has led to the development of several methods that can be broadly classified as approximate policy iteration (API) methods, where the termination criterion of conventional policy iteration is altered in order to terminate in finite time and generate an approximately optimal solution.

There are several excellent papers that consider different variants of API. However, the vast majority of these are limited to the case where the cost function is discounted, i.e.~where future costs are successively weighted less and less \citep{scherrer_approximate_2014,bertsekas_approximate_2011,stachurski_continuous_2008,santos_convergence_2004}. Though a discounted cost function may be relevant for some problems and allows for more easily determining a termination criterion, a sizeable portion of optimal control problems are better formulated as undiscounted problems (e.g.~minimum fuel/energy/time problems, or yield maximisation for chemical plants and cultivation). \citet{guo_policy_2017} introduce one API method for the undiscounted case from a reinforcement learning perspective, but this method is limited both in that the cost function must be a sum of a positive definite function of the state and a quadratically weighted function of the controls, and that the state and control cannot be arbitrarily constrained.

In this paper we will introduce a method similar to API schemes that approximates the solution to the infinite-horizon problem by instead solving a finite-horizon problem. More specifically, the method uses conventional interpolating ADP to approximate the undiscounted, infinite-horizon, non-linear, optimal control problem where the state is constrained to converge to a unique equilibrium. The primary contribution of this paper is a termination criterion that terminates at a suitable horizon without requiring the presence of a discount factor, while also allowing for (nearly) arbitrary cost, constraints, and problem dynamics --- a combination that is novel to the best of the authors knowledge. The method's sole tuning parameter allows for controlling the trade off between memory consumption, computational time, and accuracy. This allows for the method to be used without in-depth knowledge of the method. Furthermore, as the method's output is the optimal control policy (i.e.~the optimal control tabulated by the system state) subsequent on-line control can be implemented using a computationally fast interpolation operation.

The structure of this paper follows; in \secref{Problem-formulation} we will define the problem studied in this paper and the structure of the interpolating ADP method we subsequently base our presented method on. We will assume a working knowledge of ADP methods for optimal control. \citet{Sundstrom2009} gives a straightforward introduction while \citet{Bertsekas2017,Puterman1994} go into more detail. This is followed by \secref{prob-properties}, where we derive relevant properties of the studied problem. Though these properties are mostly already known, by deriving them we can both highlight some important details, as well as use a language and notation more commonly seen by control engineers as compared to existing API literature. In \secref{The-ucpadp-method} we present our method of generating an approximation of the optimal control policy, as well as highlight how existing API methods compare with our method. Finally, in \secref{Representative-examples} we use two representative examples to show the results generated by our method. For ease of reference, a list of the symbols and notation used in this paper is shown in \tabref{symbs-and-notation}.

\begin{table*}
\caption{List of used notation, symbols, and first definition.\label{tab:symbs-and-notation}}

\begin{tabular}{ccl}
DP &  & Dynamic programming\tabularnewline
ADP &  & Approximate dynamic programming\tabularnewline
API &  & Approximate policy iteration\tabularnewline
$d_{x}$ & \eqref{def_dx} & Distance between neighbouring points in $\mathcal{X}$\tabularnewline
$d_{u}$ & \eqref{def_du} & Distance between neighbouring points in $\mathcal{U}$\tabularnewline
$f_{c}$ & \eqref{base_problem_cost} & Cost function\tabularnewline
$f_{c,R}$ & \eqref{def_relaxed_cost} & Relaxed cost function\tabularnewline
$f_{d}$ & \eqref{sys_dynamics} & System dynamics function\tabularnewline
$f_{\alpha}$ & \eqref{base_problem_avg_constraint} & Average constraint function\tabularnewline
$\mathcal{F}$ & \eqref{feas_init_cond} & Set of initial conditions with feasible initial condition\tabularnewline
$\mathcal{F}'_{k}$ & \eqref{feasible-gridded-initial} & Set of feasible gridded initial conditions after $k$ samples\tabularnewline
$g$ & \eqref{base_problem_ord_constraint} & Inequality constraint function\tabularnewline
$J$ & \eqref{base_problem_cost} & Cost\tabularnewline
$J^{*}$ & \eqref{base_problem_opt_cost} & Optimal cost\tabularnewline
$J_{\text{eq}}^{*}$ & \eqref{eq_min_cost} & Optimal equilibrium cost\tabularnewline
$J_{R}$ & \eqref{relaxed_prob_cost} & Relaxed cost\tabularnewline
$J_{R}^{*}$ & \eqref{relaxed_prob_cost_opt} & Optimal relaxed cost\tabularnewline
$J_{R}^{*N}$ & \eqref{finite-horz-prob-def-cost} & Optimal relaxed $N$-horizon cost\tabularnewline
$N_{M}$ & \eqref{cond_approx_min_stationary_nm} & Finite minimum horizon\tabularnewline
$N_{M}'$ & \eqref{ucpadp-approx-horzlen} & Finite UCPADP horizon\tabularnewline
$\mathcal{S}$ & \eqref{base_problem_ord_constraint} & Set of trajectories with feasible dynamics and inequality constraints\tabularnewline
$u_{k}$ & \eqref{sys_dynamics} & Control signal at sample $k$\tabularnewline
$\bar{u}$ & \eqref{ctrl_traj} & Control trajectory\tabularnewline
$\bar{u}^{*}$ & \eqref{base_problem_opt_arg} & Optimal control trajectory\tabularnewline
$u_{\text{eq}},u_{\text{eq}}^{*}$ & \eqref{def_eq} & Optimal equilibrium control (identical by \thmref{eq_conv})\tabularnewline
$\bar{u}_{R}$ & \eqref{relaxed_prob_cost} & Relaxed control trajectory\tabularnewline
$\bar{u}_{R}^{*}$ & \eqref{relaxed_prob_traj_opt} & Optimal relaxed control trajectory\tabularnewline
$\bar{u}_{R}^{*N}$ & \eqref{finite-horz-prob-def-traj} & Optimal relaxed $N$-horizon control trajectory\tabularnewline
$\mathcal{U}$ & \eqref{def_grid_u} & Cartesian grid of sampled controls for ADP routine\tabularnewline
$\mathcal{V}_{\alpha}$ & \eqref{base_problem_avg_constraint} & Set of trajectories satisfying average equality constraint\tabularnewline
$x_{k}$ & \eqref{sys_dynamics} & System state at sample $k$\tabularnewline
$\bar{x}$ & \eqref{state_traj} & State trajectory\tabularnewline
$\bar{x}^{*}$ & \eqref{base_problem_opt_arg} & Optimal state trajectory\tabularnewline
$x_{\text{eq}}$,$x_{\text{eq}}^{*}$ & \eqref{def_eq} & Optimal equilibrium state (identical by \thmref{eq_conv})\tabularnewline
$\bar{x}_{R}$ & \eqref{relaxed_prob_cost} & Relaxed state trajectory\tabularnewline
$\bar{x}_{R}^{*}$ & \eqref{relaxed_prob_traj_opt} & Optimal relaxed state trajectory\tabularnewline
$\bar{x}_{R}^{*N}$ & \eqref{finite-horz-prob-def-traj} & Optimal relaxed $N$-horizon state trajectory\tabularnewline
$x_{k,CL}$ & \eqref{finite-horz-state} & Closed-loop state after applying a control policy $k$ times\tabularnewline
$\mathcal{X}$ & \eqref{def_grid_x} & Cartesian grid of sampled states for ADP routine\tabularnewline
$\alpha$ & \eqref{base_problem_avg_constraint} & Average constraint\tabularnewline
$\Delta_{\mu}^{k}$ & \eqref{cond-mu-test} & Control policy deviation at sample $k$\tabularnewline
$\Delta_{x}^{k}$ & \eqref{delta_x} & State deviation at sample $k$\tabularnewline
$\varepsilon_{x}$ & \eqref{cond_approx_min_stationary_eps_x} & State tolerance\tabularnewline
$\varepsilon_{\mu}$ & \eqref{cond_approx_min_stationary_eps_mu} & Control policy tolerance\tabularnewline
$\lambda$ & \eqref{def_relaxed_cost} & Relaxation parameter\tabularnewline
$\mu^{*}$ & \eqref{const_law} & Optimal stationary control law\tabularnewline
$\bar{\mu}_{R}^{*N}$ & \eqref{finite-horz-ctrl-laws} & Optimal relaxed $N$-horizon control policies\tabularnewline
\end{tabular}
\end{table*}

\section{Problem formulation\label{sec:Problem-formulation}}

Assume a dynamic system $f_{d}:\mathbb{R}^{n}\times\mathbb{R}^{m}\to\mathbb{R}^{n}$ whose associated state evolution is recursively given by
\begin{equation}
x_{k+1}=f_{d}\left(x_{k},u_{k}\right)\label{eq:sys_dynamics}
\end{equation}
for the system state $x_{k}\in\mathbb{R}^{n}$ and control input $u_{k}\in\mathbb{R}^{m}$ at samples $k\in\left[0,1,2,\dots\right]$. Define the infinite sequences 
\begin{subequations}
\begin{align}
\bar{x} & \triangleq\left[x_{0},x_{1},x_{2},\dots\right]\label{eq:state_traj}\\
\bar{u} & \triangleq\left[u_{0},u_{1},u_{2},\dots\right]\label{eq:ctrl_traj}
\end{align}
\end{subequations}
 as the \emph{state trajectory }and \emph{control trajectory} respectively. Similarly, define the finite sequences $\bar{x}^{N}\triangleq\left[x_{0},x_{1},\dots,x_{N-1}\right]$ and $\bar{u}^{N}\triangleq\left[u_{0},u_{1},\dots,u_{N-1}\right]$. In particular, for both $\bar{x}$ and $\bar{x}^{N}$ we respectively define $x_{0}$ as the \emph{initial condition}.

\subsection{The infinite-horizon problem}

Given $x_{0}$, introduce

\begin{subequations}
\label{eq:base_problem}
\begin{align}
J\left(\bar{x},\bar{u}\right) & =\lim_{N\to\infty}\frac{1}{N}\sum_{k=0}^{N-1}f_{c}\left(x_{k},u_{k}\right)\label{eq:base_problem_cost}\\
J^{*} & =\min_{\bar{x},\bar{u}}J\left(\bar{x},\bar{u}\right)\label{eq:base_problem_opt_cost}\\
\left(\bar{x}^{*},\bar{u}^{*}\right) & =\underset{\bar{x},\bar{u}}{\argmin}J\left(\bar{x},\bar{u}\right)\label{eq:base_problem_opt_arg}\\
\shortintertext{\text{subject to}}\left(\bar{x},\bar{u}\right) & \in\mathcal{S}\cap\mathcal{V}_{\alpha}\label{eq:base_problem_constraint}
\end{align}
\vspace{-2\baselineskip}
\begin{align}
\shortintertext{\text{for}}\mathcal{S} & =\left\{ \left(\bar{x},\bar{u}\right):\lim_{N\to\infty}\begin{array}{c}
g\left(x_{k},u_{k}\right)\leq0\\
x_{k+1}=f_{d}\left(x_{k},u_{k}\right)
\end{array},\,\forall k\in\left[0,N-1\right]\right\} \label{eq:base_problem_ord_constraint}\\
\mathcal{V}_{\alpha} & =\left\{ \left(\bar{x},\bar{u}\right):\lim_{N\to\infty}\frac{1}{N}\sum_{k=0}^{N-1}f_{a}\left(x_{k},u_{k}\right)=\alpha\right\} \label{eq:base_problem_avg_constraint}
\end{align}
\end{subequations}
as the problem we study in this paper. Here, we denote $f_{c}:\mathbb{R}^{n}\times\mathbb{R}^{m}\to\mathbb{R}$ the \emph{cost function}, $g:\mathbb{R}^{n}\times\mathbb{R}^{m}\to\mathbb{R}^{l}$ the \emph{inequality constraint(s)}, a scalar parameter $\alpha\in\mathbb{R}$ the \emph{average constraint}, and $f_{a}:\mathbb{R}^{n}\times\mathbb{R}^{m}\to\mathbb{R}$ the \emph{average constraint} \emph{function}. We define a \emph{feasible} \emph{trajectory} as any trajectory $\left(\bar{x},\bar{u}\right)$ that satisfies (\ref{eq:base_problem_constraint}). The set $\mathcal{S}$ gives a convenient notation for demanding that the ``textbook'' problem dynamics and inequality constraints hold, while the set $\mathcal{V}_{\alpha}$ denotes an additional average constraint.

Crucially, as none of the functions in \eqref{base_problem} are explicitly dependent on $k$, its solution satisfies the \emph{principle of optimality} (\citealt[p. 20]{Bertsekas2017}; \citealt{bellman1954theory}). \citet[p. 15]{Bertsekas2017} shows that this in turn implies that the optimal control trajectory $\bar{u}^{*}$ can equivalently be formulated as the control policy (i.e.~state-feedback) 
\begin{equation}
\bar{u}^{*}=\left[\mu_{0}^{*}\left(x_{0}\right),\mu_{1}^{*}\left(x_{1}\right),\dots\right],\label{eq:varying-state-feedback}
\end{equation}
where $\mu_{k}^{*}:\mathbb{R}^{n}\to\mathbb{R}^{m}$ are functions that are \emph{independent of the initial condition} $x_{0}$. Note that while $\bar{x}$ and $\bar{u}$ (with various sub- and super-scripts) are sequences of vectors of scalars, $\bar{\mu}$ (with various sub- and super-scripts) are instead \emph{sequences of functions}. We will refer to $\bar{\mu}^{*}$ as the \emph{optimal control policy}. 
\begin{defn}
Define $\mathcal{F}\subseteq\mathbb{R}^{n}$ as the set of initial conditions with feasible solutions, i.e.
\begin{equation}
\mathcal{F}\triangleq\left\{ x_{0}:\exists\left(\bar{x},\bar{u}\right)\in\mathcal{S}\cap\mathcal{V}_{\alpha}\right\} .\label{eq:feas_init_cond}
\end{equation}
\end{defn}

\begin{assm}For the remainder this paper we assume:\label{assm:base-assumptions}
\begin{enumerate}[label=\bfseries{A.\arabic{enumi}}]
\item $f_{c}$, $f_{d}$, $g$, and $f_{a}$ are continuous and bounded. \label{enu:base-assumptions_bounded}
\item The optimal solution $\left(\bar{x}^{*},\bar{u}^{*}\right)$ associated with $x_{0}$ is unique.\label{enu:unique_sol}
\item The optimal control policy associated with \eqref{base_problem} exists, and can be expressed as
\begin{equation}
\bar{u}^{*}=\left[\mu^{*}\left(x_{0}\right),\mu^{*}\left(x_{1}\right),\dots\right],\label{eq:const_law}
\end{equation}
i.e.~it is not only independent of the initial condition $x_{0}$, but also independent of the sample index $k$. We will refer to this as a \emph{stationary }control policy \citep{BERTSEKAS1979607}.\label{enu:base-assumptions-stationary}
\item $\mathcal{F}$ is nonempty, $\lim_{k\to\infty}\left(x_{k}^{*}\right)$ exists and is independent of $x_{0}$ for all $x_{0}\in\mathcal{F}$, and $x_{k}^{*}$ is asymptotically stable in the sense of Lyapunov for $x_{0}$ near $\lim_{k\to\infty}\left(x_{k}^{*}\right)$.\label{enu:base-assumptions_feas_unique_eq}
\end{enumerate}
\end{assm}

Note that \enuref{base-assumptions_bounded} implies that $J\left(\bar{x},\bar{u}\right)$ is finite for any feasible trajectory, and by \enuref{base-assumptions_feas_unique_eq} we can furthermore view $J^{*}$ as the \emph{average (mean) cost}.
\begin{defn}
Assuming \enuref{base-assumptions_feas_unique_eq} holds, define
\begin{equation}
\left(x_{\text{eq}},u_{\text{eq}}\right)\triangleq\lim_{k\to\infty}\left(x_{k}^{*},u_{k}^{*}\right)\label{eq:def_eq}
\end{equation}
as the problem's equilibrium point.
\end{defn}
Note that \enuref{unique_sol}, \enuref{base-assumptions-stationary}, and \enuref{base-assumptions_feas_unique_eq} may be difficult to determine a priori for a given problem. We will briefly discuss the possible effects of them not holding in \secref{The-ucpadp-method}.

\subsection{Interpolating ADP}

The method we introduce in this paper uses a conventional interpolating ADP scheme, and we will here use the standard method of gridding $x$ and $u$ into finite Cartesian sets. We define 
\begin{subequations}
\begin{align}
d_{x} & \in\mathbb{R}^{n}\label{eq:def_dx}\\
d_{u} & \in\mathbb{R}^{m}\label{eq:def_du}
\end{align}
\end{subequations}
 as the distance between neighbouring grid points for each dimension of the states and controls respectively. We also define

\begin{subequations}
\begin{align}
\mathcal{X} & \subset\mathbb{R}^{n}\label{eq:def_grid_x}\\
\mathcal{U} & \subset\mathbb{R}^{m}\label{eq:def_grid_u}
\end{align}
\end{subequations}
 as the discrete set of state and control grid points resolved by ADP respectively, separated by $d_{x}$ and $d_{u}$ respectively and bounded by the region(s) where $g(x,u)\leq0$. We then use conventional multilinear interpolation to approximate the cost $J$ and optimal control policy $\mu$ for the real-valued states that do not lie in the discrete set $\mathcal{X}$. For example, assuming $x\in\mathbb{R}^{2}$, $u\in\mathbb{R}^{1}$, and $g(x,u)=|x|_{1}\leq1\land|u|\leq1$, choosing the very coarse (but illustrative) $d_{x}=$$\left[2,2\right]^{T}$ and $d_{u}=0.5$ gives the sets 
\begin{subequations}
\begin{align}
\mathcal{X} & =\left\{ \left[\begin{array}{c}
-1\\
-1
\end{array}\right],\left[\begin{array}{c}
-1\\
1
\end{array}\right]\left[\begin{array}{c}
1\\
-1
\end{array}\right],\left[\begin{array}{c}
1\\
1
\end{array}\right]\right\} \\
\mathcal{U} & =\left\{ -1,-0.5,0,0.5,1\right\} .
\end{align}
\end{subequations}

\section{Infinite-horizon, average-constrained problem properties\label{sec:prob-properties}}

In this section we introduce properties of the undiscounted, infinite-horizon, average-constrained problem that will later be utilised by the method we introduce in \secref{The-ucpadp-method}.

\subsection{Solution convergence}
\begin{defn}
For $x\in\mathbb{R}^{n}$, $u\in\mathbb{R}^{m}$, using the same functions as in \eqref{base_problem}, define
\begin{subequations}
\label{eq:eq_min}
\begin{align}
J_{\text{eq}}^{*} & =\min_{x,u}f_{c}\left(x,u\right)\label{eq:eq_min_cost}\\
\left(x_{\text{eq}}^{*},u_{\text{eq}}^{*}\right) & =\underset{x,u}{\argmin}f_{c}\left(x,u\right)\label{eq:eq_min_arg}\\
\shortintertext{\text{subject to}}x & =f_{d}\left(x,u\right)\label{eq:eq_min_stationary_state}\\
g\left(x,u\right) & \leq0\\
f_{a}\left(x,u\right) & =\alpha\\
\forall x_{0}\in\mathcal{F},\exists\left(\bar{x},\bar{u}\right)\text{ s.t.} & \lim_{k\to\infty}\left(x_{k},u_{k}\right)=\left(x,u\right)\label{eq:eq_min_reachability}
\end{align}
\end{subequations}
as the \emph{optimal reachable equilibrium operating point} $\left(x_{\text{eq}}^{*},u_{\text{eq}}^{*}\right)$. (Note that we have identical states on both the left- and right-hand side of (\ref{eq:eq_min_stationary_state}), i.e.~an equilibrium state.) We can view this as the unique stationary point of the system with lowest cost that we can reach for any initial condition in the feasible set $\mathcal{F}$.
\end{defn}
\begin{thm}
\label{thm:eq_conv} Given \enuref{base-assumptions_bounded} and \enuref{base-assumptions_feas_unique_eq},
\begin{align}
J^{*} & =J_{\text{eq}}^{*}\\
\left(x_{\text{eq}},u_{\text{eq}}\right) & =\left(x_{\text{eq}}^{*},u_{\text{eq}}^{*}\right),
\end{align}
i.e.~the equilibrium we reach will be optimal in the sense of (\ref{eq:eq_min}).
\end{thm}
\begin{IEEEproof}
For $0\leq i<j$, define 
\begin{equation}
J_{i\to j}\left(\bar{x},\bar{u}\right)\triangleq\sum_{k=i}^{j}f_{c}\left(x_{k},u_{k}\right).
\end{equation}
We can then formulate \eqref{base_problem_opt_cost} as
\begin{equation}
J^{*}=\min_{\bar{x},\bar{u}}\lim_{N\to\infty}\frac{1}{N}J_{0\to i-1}\left(\bar{x},\bar{u}\right)+\frac{1}{N}J_{i\to N-1}\left(\bar{x},\bar{u}\right).
\end{equation}
As $N\to\infty$, we are guaranteed that $\frac{1}{N}J_{0\to i-1}=0$ for any fixed $i>0$ per our assumption that $f_{c}$ is bounded. This implies that $J^{*}$ is only dependent on $J_{i\to N-1}$. By \enuref{base-assumptions_feas_unique_eq}, we can make $\left(x_{i}^{*},u_{i}^{*}\right)$ arbitrarily close to $\left(x_{\text{eq}},u_{\text{eq}}\right)$ for sufficiently large $i$.

Suppose that 
\begin{equation}
\left(x_{\text{eq}},u_{\text{eq}}\right)\neq\left(x_{\text{eq}}^{*},u_{\text{eq}}^{*}\right).\label{eq:eq_conv_contradiction}
\end{equation}
By \enuref{base-assumptions_feas_unique_eq} $\left(x_{\text{eq}},u_{\text{eq}}\right)$ is unique, implying that $J^{*}>J_{\text{eq}}^{*}$. However, by \eqref{eq_min_reachability} there exists trajectories $\bar{x}'$ and $\bar{u}'$ such that $\lim_{k\to\infty}\left(x'_{k},u'_{k}\right)=\left(x_{\text{eq}}^{*},u_{\text{eq}}^{*}\right)$, with corresponding cost $J'<J^{*}$, contradicting \eqref{eq_conv_contradiction}.

For an alternate view of the same proof, see \citet[p. 298]{Bertsekas2012}.
\end{IEEEproof}
By \thmref{eq_conv} we can intuitively view the infinite-horizon problem's solution as ignoring any (finite) costs during the transient phase and driving the state to the reachable stationary point with lowest cost. This is a special case of the \emph{turnpike property }\citep{TRELAT201581,Zaslavski2014}, which states that the solution to problems with a sufficiently long (finite) horizon tends to display transient dynamic initial and terminal phases, with a middle stationary phase that is independent of the initial and terminal conditions. Of course, the infinite-horizon problem does not have a terminal phase, and we can thus view the solution to our problem \eqref{base_problem} as consisting of an initial transient followed by stationary operation at the optimal reachable equilibrium point. 

From a notation perspective, by \thmref{eq_conv} we do not need to make the distinction between $x_{\text{eq}}$ and $x_{\text{eq}}^{*}$. For consistency, we will use $x_{\text{eq}}^{*}$ from here on out.

\subsection{Average-constraint relaxation\label{subsec:Average-constraint-relaxation}}
\begin{defn}
For a fixed, bounded, scalar relaxation parameter $\lambda\in\mathbb{R}$, define the \emph{relaxed cost} as
\begin{equation}
f_{c,R}\left(x,u\right)\triangleq f_{c}\left(x,u\right)+\lambda f_{a}\left(x,u\right).\label{eq:def_relaxed_cost}
\end{equation}
Now we can introduce the \emph{relaxed problem} as 
\begin{subequations}
\label{eq:relaxed_prob}
\begin{align}
J_{R}\left(\bar{x}_{R},\bar{u}_{R}\right) & =\lim_{N\to\infty}\frac{1}{N}\sum_{k=0}^{N-1}f_{c,R}\left(x_{k},u_{k}\right)\label{eq:relaxed_prob_cost}\\
J_{R}^{*} & =\min_{\bar{x}_{R},\bar{u}_{R}}J_{R}\left(\bar{x}_{R},\bar{u}_{R}\right)\label{eq:relaxed_prob_cost_opt}\\
\left(\bar{x}_{R}^{*},\bar{u}_{R}^{*}\right) & =\underset{\bar{x}_{R},\bar{u}_{R}}{\argmin}J_{R}\left(\bar{x}_{R},\bar{u}_{R}\right)\label{eq:relaxed_prob_traj_opt}\\
\shortintertext{\text{subject to}}\left(\bar{x}_{R},\bar{u}_{R}\right) & \in\mathcal{S},\label{eq:relaxed_prob_constraint}
\end{align}
\end{subequations}
where we view $J_{R}\left(\bar{x}_{R},\bar{u}_{R}\right)$ as the \emph{relaxed representation }of $J\left(\bar{x},\bar{u}\right)$, and $J_{R}^{*}$ and $\left(\bar{x}_{R}^{*},\bar{u}_{R}^{*}\right)$ as the \emph{optimal relaxed cost} and \emph{optimal relaxed trajectories} respectively. Note that $\left(\bar{x}_{R},\bar{u}_{R}\right)$, and therefore also $\left(\bar{x}_{R}^{*},\bar{u}_{R}^{*}\right)$, are not formally constrained to lie in $\mathcal{V}_{\alpha}$. 

For clarity, we will use the notation $\bar{x}_{R}$ and $\bar{u}_{R}$ when referring to trajectories associated with the relaxed problem. We will for ease of notation assume that $\left(\bar{x}_{R}^{*},\bar{u}_{R}^{*}\right)$ is unique (much as \eqref{base_problem_opt_arg}), though we can in principle use DP (and in turn the method to be presented) to solve problems with non-unique solutions.
\end{defn}
\begin{lem}
For a given $\alpha$, assume for some $\lambda$ we have $\left(\bar{x}_{R}^{*},\bar{u}_{R}^{*}\right)\in\mathcal{V}_{\alpha}$. Then $\left(\bar{x}_{R}^{*},\bar{u}_{R}^{*}\right)=\left(\bar{x}^{*},\bar{u}^{*}\right)$.\label{lem:relaxation}
\end{lem}
\begin{IEEEproof}
For convenience, introduce $\zeta^{*}\triangleq\left(\bar{x}^{*},\bar{u}^{*}\right)$, $\zeta_{R}^{*}\triangleq\left(\bar{x}_{R}^{*},\bar{u}_{R}^{*}\right)$, $\zeta\triangleq\left(\bar{x},\bar{u}\right)$, $\zeta_{R}\triangleq\left(\bar{x}_{R},\bar{u}_{R}\right)$, and 
\begin{equation}
h\left(\zeta\right)\triangleq\lim_{N\to\infty}\frac{1}{N}\sum_{k=0}^{N-1}f_{a}\left(x_{k},u_{k}\right)-\alpha.
\end{equation}
Note that $h\left(\zeta\right)=0\Leftrightarrow\zeta\in\mathcal{V}_{\alpha}.$

The weak duality theorem \citep{andreasson-intro-continuous-optimization-2016} ensures that 
\begin{align}
J\left(\zeta^{*}\right) & \geq J\left(\zeta_{R}^{*}\right)+\lambda h\left(\zeta_{R}^{*}\right)=J_{R}\left(\zeta_{R}^{*}\right)-\lambda\alpha.\label{eq:relaxed-weak-duality-relaxed}
\end{align}

In \eqref{relaxed-weak-duality-relaxed}, by our assumption $\zeta_{R}^{*}\in\mathcal{V}_{\alpha}$ we are ensured that $h\left(\zeta_{R}^{*}\right)=0$, giving
\begin{equation}
J\left(\zeta^{*}\right)\geq J\left(\zeta_{R}^{*}\right)=J_{R}\left(\zeta_{R}^{*}\right)-\lambda\alpha.\label{eq:relaxed-weak-duality-relaxed-2}
\end{equation}
As $\zeta_{R}^{*}\in\mathcal{S}\cap\mathcal{V}_{\alpha}$, $\zeta_{R}^{*}$ also minimises \eqref{base_problem}, allowing us to replace the inequality in \eqref{relaxed-weak-duality-relaxed-2} with strict equality. By \enuref{unique_sol} $\zeta^{*}$ and $\zeta_{R}^{*}$ are unique, ensuring that that $\zeta^{*}=\zeta_{R}^{*}$. Finally, as $\zeta_{R}^{*}$ is independent of constant terms we have that 
\begin{equation}
\zeta^{*}=\zeta_{R}^{*}=\underset{\zeta}{\argmin}J_{R}\left(\zeta\right).
\end{equation}
\end{IEEEproof}
\begin{thm}
\label{thm:relaxed}Given \eqref{base_problem_opt_arg} and its relaxed counterpart \eqref{relaxed_prob_traj_opt},
\begin{enumerate}[label=\bfseries{R.\arabic{enumi}}]
\item If \eqref{relaxed_prob_traj_opt} is infeasible (i.e.~a solution does not exist), then \eqref{base_problem_opt_arg} is also infeasible (i.e.~\enuref{base-assumptions_feas_unique_eq} is violated).\label{enu:relaxed-infeas}
\item For a given $\lambda$ and feasible \eqref{relaxed_prob_traj_opt}, there exists an $\alpha$ where\label{enu:relaxed-optimal}
\end{enumerate}
\begin{equation}
\left(\bar{x}_{R}^{*},\bar{u}_{R}^{*}\right)=\left(\bar{x}^{*},\bar{u}^{*}\right).\label{eq:relaxed_solution_equal}
\end{equation}

\end{thm}
\begin{IEEEproof}
\enuref{relaxed-infeas}: Trivial, as $\mathcal{S}\supseteq\mathcal{S}\cap\mathcal{V}_{\alpha}$ .
\end{IEEEproof}

\begin{IEEEproof}
\enuref{relaxed-optimal}: As $\lambda$ is given and \eqref{relaxed_prob_traj_opt} is feasible, we can thus find $\left(\bar{x}_{R}^{*},\bar{u}_{R}^{*}\right)$. Let us now define
\begin{equation}
\alpha'\triangleq\lim_{N\to\infty}\frac{1}{N}\sum_{k=0}^{N-1}f_{a}\left(\bar{x}_{R,k}^{*}\bar{u}_{R,k}^{*}\right).\label{eq:relaxed_alpha_def}
\end{equation}
For $\alpha=\alpha'$ we (by construction) have $(\bar{x}_{R}^{*},\bar{u}_{R}^{*})\in\mathcal{S}\cap\mathcal{V}_{\alpha}$, trivially satisfying the requirements of \lemref{relaxation}.
\end{IEEEproof}
In essence, for a given $\lambda$ \enuref{relaxed-optimal} ensures us that $(\bar{x}_{R}^{*},\bar{u}_{R}^{*})=(\bar{x}^{*},\bar{u}^{*})$ for some value of $\alpha$. We can intuitively view $\lambda$ as a tuning parameter, where different values of $\lambda$ are associated with different solutions, each of which (trivially) have an associated average that we can compute by means of \eqref{relaxed_alpha_def}.

Using the relaxed problem formulation allows us to avoid the explicit average constraint \eqref{base_problem_avg_constraint}, which is primarily of use in the sense that the problem becomes more numerically tractable. At its core, the method we will introduce in this paper approximates the solution to \eqref{base_problem} by instead solving a finite-horizon problem of sufficient length. One naive method of satisfying the average constraint would then be to introduce an additional state variable that stores the accumulated average, i.e. $z_{N}=\sum_{k=0}^{N-1}f_{\alpha}\left(x_{k},u_{k}\right)$. We could then add an equality constraint demanding $z_{n}/N=\alpha$. However, this is computationally demanding (as we need to introduce an additional state variable, which DP schemes scale poorly with) and introduces a bias in the achieved average (as the average $z_{n}/N$ is taken over both the initial transient and the stationary phase, we therefore only achieve the desired average as $N\to\infty$). Using the relaxed formulation thus avoids these issues entirely.

\subsection{Convergence of finite-horizon problem }

We will in this section introduce notation for the finite-horizon problem, which will then be used for constructing the method presented in this paper.
\begin{defn}
For a given finite horizon $N$, bounded $\lambda$, and initial condition $x_{0}$, define
\begin{subequations}
\label{eq:finite-horz-prob-def}
\begin{align}
J_{R}^{*N} & =\min_{\bar{x}_{R}^{N},\bar{u}_{R}^{N}}\frac{1}{N}\sum_{k=0}^{N-1}f_{c,R}\left(x_{k},u_{k}\right)\label{eq:finite-horz-prob-def-cost}\\
\left(\bar{x}_{R}^{*N},\bar{u}_{R}^{*N}\right) & =\underset{\bar{x}_{R}^{N},\bar{u}_{R}^{N}}{\argmin}\frac{1}{N}\sum_{k=0}^{N-1}f_{c,R}\left(x_{k},u_{k}\right)\label{eq:finite-horz-prob-def-traj}\\
\shortintertext{\text{subject to}}\left(\bar{x}_{R}^{N},\bar{u}_{R}^{N}\right) & \in\mathcal{S}
\end{align}
\end{subequations}
as the \emph{N-horizon relaxed problem} with average cost $J_{R}^{*N}$ and associated (finite-length) state and control trajectories $\left(\bar{x}_{R}^{*N},\bar{u}_{R}^{*N}\right)$. Furthermore, define
\begin{equation}
\bar{\mu}_{R}^{*N}=\left[\mu_{R,0}^{*N},\mu_{R,1}^{*N},\dots,\mu_{R,N-1}^{*N}\right],\label{eq:finite-horz-ctrl-laws}
\end{equation}
where $\mu_{R,k}^{*N}:\mathbb{R}^{n}\to\mathbb{R}^{m}$ is the $k$'th state-feedback control policy, as the $N$-horizon sequence of control policies associated with (\ref{eq:finite-horz-prob-def}).
\end{defn}

\begin{defn}
Define 
\begin{equation}
x_{k,CL}\left(\mu,x_{0}\right)\label{eq:finite-horz-state}
\end{equation}
as the (not necessarily optimal) $k$'th closed-loop state given by repeatedly applying a (sample-independent) control policy $\mu$ $k$ times from an initial state $x_{0}$, e.g.
\begin{subequations}
\begin{align*}
x_{0,CL}\left(\mu,x_{0}\right) & \triangleq x_{0}\\
x_{1,CL}\left(\mu,x_{0}\right) & \triangleq f_{d}\left(x_{0},\mu\left(x_{0}\right)\right)\\
x_{2,CL}\left(\mu,x_{0}\right) & \triangleq f_{d}\left(x_{1,CL}\left(\mu,x_{0}\right),\mu\left(x_{1,CL}\left(\mu,x_{0}\right)\right)\right).
\end{align*}
\end{subequations}
 Note that the method of generating $x_{k,CL}$ is very similar to the forward-calculation stage of ADP, and differs only in that the control policy is kept constant.
\end{defn}

\begin{defn}
For a given control policy $\mu$, define 
\begin{equation}
\mathcal{F}'_{k}(\mu)\triangleq\left\{ x_{0}\in\mathcal{X}:g(x_{k',CL},\mu(x_{k',CL}))\leq0\,\forall k'\in[0,k]\right\} .\label{eq:feasible-gridded-initial}
\end{equation}
We can thus view $\mathcal{F}'_{k}(\mu)$ as the set of initial conditions in $\mathcal{X}$ that satisfies the problem constraints and dynamics (the latter trivially, as we use $\mu$ to apply a control and give the next state) after applying the control policy $\mu$ $k$ times.
\end{defn}

\begin{defn}
For $k>0$, introduce the maximum control policy deviation $\Delta_{\mu}^{k}\in\mathbb{R}^{m}$ as\label{def:delta_mu}
\begin{equation}
\left[\Delta_{\mu}^{k}\right]_{i}\triangleq\max_{\begin{array}{c}
x\in\mathcal{F}'_{\left\lceil k/2\right\rceil }(\mu_{R,0}^{*k})\\
k'\in\left[0,\left\lceil k/2\right\rceil \right]
\end{array}}\left|\left[\mu_{R,0}^{*k}\left(x\right)-\mu_{R,k'}^{*k}\left(x\right)\right]_{i}\right|,\label{eq:cond-mu-test}
\end{equation}
where the notation $\left[a\right]_{i}$ refers to the $i$'th element of a vector $a$ and $\left\lceil \dots\right\rceil $ refers to the ceiling function. We can view $\Delta_{\mu}^{k}$ as indicating the convergence of $\mu_{R,0}^{*k}$ to $\mu^{*}$, evaluated at the gridded state points $\mathcal{X}$ whose associated state evolution remains feasible after $k/2$ iterations.
\end{defn}

\begin{defn}
Introduce the maximum state deviation $\Delta_{x}^{k}\in\mathbb{R}^{n}$ as\label{def:delta_x}
\begin{align}
\left[\Delta_{x}^{k}\right]_{i} & \triangleq\max_{x\in\mathcal{F}'_{k\left\lceil k/2\right\rceil }(\mu_{R,0}^{*k})}\Biggl|\Biggl[x_{\left\lceil k/2\right\rceil ,CL}\left(\mu_{R,0}^{*k},x\right)-\nonumber \\
 & \sum_{x'\in\mathcal{F}'_{\left\lceil k/2\right\rceil }(\mu_{R,0}^{*k})}x_{\left\lceil k/2\right\rceil ,CL}\left(\mu_{R,0}^{*k},x'\right)\frac{1}{|\mathcal{F}'_{\left\lceil k/2\right\rceil }|}\Biggr]_{i}\Biggr|.\label{eq:delta_x}
\end{align}
Note that the notationally heavy second line of \eqref{delta_x} is equivalent to the mean feasible state after $\left\lceil k/2\right\rceil $ iterations. Similarly to \defref{delta_mu}, we can thus view $\Delta_{x}^{k}$ as indicating the convergence of $\left[x_{0,CL},x_{1,CL},\dots,x_{\left\lceil k/2\right\rceil ,CL}\right]$ to $\bar{x}^{*}$, evaluated at the points where $x_{\left\lceil k/2\right\rceil ,CL}$ remains feasible.
\end{defn}
Trivially, using \defref{delta_mu} and \defref{delta_x} gives:
\begin{prop}
By \enuref{base-assumptions-stationary} $\lim_{k\to\infty}\Delta_{\mu}^{k}=0$, and by \enuref{base-assumptions_feas_unique_eq} $\lim_{k\to\infty}\Delta_{x}^{k}=0$.\label{prop:conv_ucpadp_eps}
\end{prop}
\begin{defn}
Given a \emph{control policy tolerance} $\varepsilon_{\mu}\in\mathbb{R}^{m}$ and \emph{state convergence tolerance} $\varepsilon_{x}\in\mathbb{R}^{n}$, define

\begin{subequations}
\label{eq:cond_approx_min_stationary}
\begin{align}
N_{M} & \triangleq\min_{k}k\label{eq:cond_approx_min_stationary_nm}\\
\shortintertext{\text{such that}}\left[\Delta_{\mu}^{k}\right]_{i} & <\left[\varepsilon_{\mu}\right]_{i}\forall i\in\left[1,m\right]\label{eq:cond_approx_min_stationary_eps_mu}\\
\left[\Delta_{x}^{k}\right]_{i} & <\left[\varepsilon_{x}\right]_{i}\forall i\in\left[1,n\right],\label{eq:cond_approx_min_stationary_eps_x}
\end{align}
\end{subequations}
as the \emph{minimum horizon}. Proposition~\propref{conv_ucpadp_eps} ensures us that that for any $\varepsilon_{\mu}$ and $\varepsilon_{x}$ there exists an associated finite horizon $N_{M}$, which we view as the shortest finite-horizon approximation of the infinite-horizon problem. 
\end{defn}

\section{The UCPADP method\label{sec:The-ucpadp-method}}

In this section we introduce the primary contribution of this paper: \emph{Undiscounted Control Policy generation by Approximate Dynamic Programming} (UCPADP), a method that generates an approximation of $\mu^{*}$. At its core, in UCPADP we generate an approximation of the optimal control policy by iteratively testing successively larger horizons until the termination criteria (\ref{eq:cond_approx_min_stationary}) are satisfied. For computational efficiency reasons we will return to, UCPADP will approximate the control policy as 
\begin{subequations}
\begin{align}
\mu^{*} & \approx\mu_{R,0}^{*N_{M}'}\label{eq:ucpadp-approx}\\
\text{where }N_{M} & \leq N_{M}'\leq2N_{M},\label{eq:ucpadp-approx-horzlen}
\end{align}
\end{subequations}
i.e.~the generated horizon will lie in a range between $N_{M}$ and $2N_{M}$.

We can at this stage highlight one of the more significant differences between UCPADP and conventional API: the choice of termination conditions. Conventional API generates improved control policies analogous to $\mu_{R,1}^{*1},\mu_{R,2}^{*2},\mu_{R,3}^{*3},\dots$ with an associated cost $J_{R}^{*1},J_{R}^{*2},J_{R}^{*3},\dots$, and eventually terminates when the difference between either successive policies or cost is below a given threshold, for instance as in \citet[CPI, PSDP]{santos_convergence_2004,stachurski_continuous_2008}. This is similar to the test performed in \eqref{cond_approx_min_stationary_eps_mu}, which requires the control policy to be near-stationary. However, in conventional API the termination tolerance (analogous to $\varepsilon_{\mu}$) is sized based on the discount factor, and depending on the specific method chosen the tolerance is either undefined or tends towards zero when the discount factor tends towards one (i.e.~becomes the undiscounted case we study here). \citet{scherrer_approximate_2014,bertsekas_approximate_2011} review other methods that do not terminate based on the change in the control policy, but instead use some other termination criterion. However, these methods also assume a discounted problem formulation. \citet{guo_policy_2017} is one example of a method that considers the undiscounted case, however their method imposes fairly significant limits on the class of cost and constraint functions (as discussed previously).

The state convergence condition \eqref{cond_approx_min_stationary_eps_x} is to the best of our knowledge novel, and serves a crucial purpose in that it demands the horizon be long enough for \emph{all gridded feasible initial conditions to converge to a region near the equilibrium}. Recall that by \enuref{base-assumptions_feas_unique_eq} $x_{k}^{*}$ (the true optimal state trajectory) is stable in the sense of Lyapunov for initial conditions near the equilibrium, and in concert with \thmref{eq_conv} we are thus ensured that an initial condition near the equilibrium will also remain in its vicinity. As we apply test \eqref{cond_approx_min_stationary_eps_x} to all feasible elements in $\mathcal{X}$, at least one initial condition $x_{0}\in\mathcal{X}$ will therefore start and then \emph{remain in the nearby vicinity of the equilibrium}. Ultimately, by combining \eqref{cond_approx_min_stationary_eps_mu} and \eqref{cond_approx_min_stationary_eps_x} we are ensured that $\mu_{R,0}^{*N_{M}'}$ is nearly constant during the interval needed for all feasible gridded points in $\mathcal{X}$ to reach the vicinity of the equilibrium.

\begin{figure}
\begin{centering}
\includegraphics[width=0.8\columnwidth]{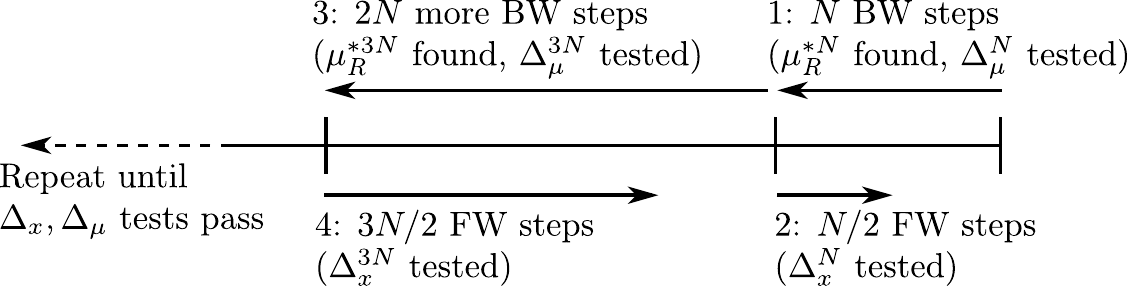}
\par\end{centering}
\caption{UCPADP steps, successively switching between generating more accurate control policies (backward calculation, steps 1,3,$\dots$), and evaluating whether the control policy is constant over the time needed for the state evolution to converge (forward calculation, steps 2,4,$\dots$).\label{fig:ucpadp-phases}}
\end{figure}

In UCPADP, we determine $\mu_{R,0}^{*N_{M}'}$ numerically efficiently in a manner similar to API implemented with ADP. We do this using a nested scheme that repeatedly switches between backward-calculation phases (successively generating control policies with longer associated horizons) and forward-calculation phases (applying tests \eqref{cond_approx_min_stationary_eps_mu} and \eqref{cond_approx_min_stationary_eps_x}, and eventually terminating when both tests pass). A description of the phases in UCPADP follows, see \figref{ucpadp-phases} for an illustration. For now, assume $\varepsilon_{\mu}$ and $\varepsilon_{x}$ are given (fixed) vectors.

First, we arbitrarily choose a small initial horizon $N$ and perform $N$ backward-calculation iterations, giving us (among other data) $\mu_{R}^{*N}$. We can then perform test \eqref{cond_approx_min_stationary_eps_mu} and, by performing $N/2$ forward-calculation steps, test \eqref{cond_approx_min_stationary_eps_x}. If both tests pass we terminate and return $\mu_{R,0}^{*N}$ as our approximation of $\mu^{*}$. Conversely, if either of these tests fail by \enuref{base-assumptions_feas_unique_eq} we are ensured that increasing the horizon sufficiently will give a control policy that satisfies the tests. In UCPADP we chose to proceed by increasing the horizon to $3N$. Fortunately, in our DP scheme we can compute $\mu_{R}^{*3N}$ using only $2N$ additional backward-calculation iterations by resuming the backward-calculation from $\mu_{R}^{*N}$. This is possible as each successive backward-calculation step is independent of the total horizon. After generating $\mu_{R}^{*3N}$ we can now again test \eqref{cond_approx_min_stationary_eps_mu} and \eqref{cond_approx_min_stationary_eps_x}. Should both tests pass we can return $\mu_{R,0}^{*3N}$ as our approximation of $\mu^{*}$, and otherwise recursively repeat this procedure of doubling the number of back-calculation steps until the tests pass (i.e.~generating and testing horizons $N,3N,9N,27N,\dots$). A pseudocode implementation of the UCPADP method is listed in Algorithm~\ref{alg:ucpadp}.

Up to this point we have assumed that the problem solution is unique (\enuref{unique_sol}), converges to a stationary control policy (\enuref{base-assumptions-stationary}), and all states converge to a unique equilibrium (\enuref{base-assumptions_feas_unique_eq}). Let us now briefly consider the case where we do not know if these assumptions hold beforehand. Beginning with \enuref{unique_sol}, recall that we can determine whether or not this assumption holds during the backward-calculation phase by checking if the minimum cost is unique, and in the case of a non-unique cost we can resolve this by simply returning one arbitrarily selected optimal solution. Let us now focus on the case where \enuref{base-assumptions-stationary} and \enuref{base-assumptions_feas_unique_eq} are unverified. Applying the UCPADP method gives one of two possible outcomes: UCPADP either never terminates (i.e.~\eqref{cond_approx_min_stationary_eps_mu} and \eqref{cond_approx_min_stationary_eps_x} never pass), or it terminates after a finite number of back-calculation iterations. If UCPADP never terminates, then one possible cause is that \enuref{base-assumptions-stationary} and/or \enuref{base-assumptions_feas_unique_eq} do not hold (i.e.~the termination criteria \eqref{cond_approx_min_stationary_eps_mu} and \eqref{cond_approx_min_stationary_eps_x} correctly detected a non-stationary control policy and/or detected that the system states do not converge to a single equilibrium). Alternatively, it is possible that the problem's discretisation and/or tolerances were poorly chosen. Regardless, should UCPADP never terminate it is clear that no valid solution could be generated. If UCPADP \emph{does} terminate, we are assured that either: \emph{(i)} \enuref{base-assumptions-stationary} and \enuref{base-assumptions_feas_unique_eq} do hold and a near-optimal control policy is generated, or \emph{(ii)} the problem is maliciously nonlinear and \enuref{base-assumptions-stationary} and/or \enuref{base-assumptions_feas_unique_eq} do not hold (which went undetected by \eqref{cond_approx_min_stationary_eps_mu} and \eqref{cond_approx_min_stationary_eps_x}), ultimately giving a control policy without any clear optimality guarantees. As the class of problems we can attempt to solve with UCPADP covers general non-linear systems it is not surprising that there exist pathological problems that lead UCPADP (and ADP in general) to generate erroneous solutions. Ultimately it is up to the user of UCPADP to determine whether or not the studied problem is of a class that satisfies the (arguably mild) assumptions \enuref{base-assumptions-stationary} and \enuref{base-assumptions_feas_unique_eq}.

In Algorithm~\ref{alg:ucpadp}, we extend the notion of termination used thus far by adding a parameter $N_{\text{max}}$ that allows for configuring a maximum horizon that terminates UCPADP if $N>N_{\text{max}}$. This acts as a safety and guarantees that UCPADP terminates after a finite number of iterations. In the event that this limit triggers UCPADP to terminate we can conclude that either the minimum horizon is larger than $N_{\text{max}}$, that \enuref{base-assumptions-stationary} and/or \enuref{base-assumptions_feas_unique_eq} do not hold, or the discretisation and/or tolerances were poorly chosen. Of course, should this happen then we can not say anything about the stability (let alone the optimality) of the returned control policy.

\begin{algorithm}
\captionsetup{labelfont={bf}, labelsep=period}\caption{Pseudocode UCPADP algorithm. Here, $\mathrm{DP_{\text{1-back}}}$ and $\mathrm{DP_{\text{1-fw}}}$ are the one-step backward and forward ADP operations. $\mathcal{X}$ is the set of initial conditions tested in the ADP method. $N_{\text{init}}$ is the initial problem horizon. $C_{N}$ is the cost-to-go after $N$ iterations. Note here that a reverse notation is used for the calculated control policy; $\mu_{1}$ corresponds to the state-feedback control policy from the first back-calculation step (i.e.~$\mu_{N-1}^{*N}$) while $\mu_{N}$ corresponds to the last (i.e.~$\mu_{0}^{*N}$). We can view the index $k$ as counting the number of back-calculation steps performed. Note the abuse of notation on line~\ref{lst:line_epstest} that indicates the $\Delta_{x}^{N}$ and $\Delta_{\mu}^{N}$ tests respectively.\label{alg:ucpadp}}
\begin{algorithmic}[1]
\small
\Function{UCPADP}{$\mathcal{X}$, $N_\text{max}$, $N_\text{init}$}
	\State $N_b \gets N_\text{init}$  	\Comment{Batch back-calculation steps}
	\State $N \gets 0$					\Comment{Cumulative back-calculation steps}
	\State $C_0 \gets 0$				  \Comment{Set initial cumulative cost to zero}
	\Repeat
		\For{$N \gets N,N+N_b$}
			\State $\mu_{N+1}, C_{N+1} \gets \mathrm{DP_\text{1-back}}(C_N)$
		\EndFor

		\State{$X_{CL} \gets \mathcal{X}$}
		\For{$i \gets 1, \lceil N/2 \rceil $}
			\State{$X_{CL} \gets \mathrm{DP_\text{1-fw}}(X_{CL}, \mu_N)$}
		\EndFor
		\State{$N_\mathrm{b} \gets 2 \cdot N_\mathrm{b}$}	\Comment{Raise $N_b$ by doubling}

	\Until{$N > N_\mathrm{max}$ or $(|X_{CL} - \mathrm{mean}(X_{CL} )| < \varepsilon_x$ and \label{lst:line_epstest} $|\mu_N - \mu_k| < \varepsilon_u \ \forall k \in [ \lceil N/2 \rceil ,N])$}
\State \textbf{return} $\mu_N$, $X_{CL}$, $N$
\EndFunction 
\end{algorithmic}
\end{algorithm}

Tests \eqref{cond_approx_min_stationary_eps_mu} and \eqref{cond_approx_min_stationary_eps_x} are straightforward to compute exhaustively, as the initial conditions $x_{0}$ come from the discrete set $\mathcal{X}$. Furthermore, in UCPADP we have chosen to double the number of additional back-calculation steps to perform between each test evaluation. This attempts to balance the time spent on backward-calculation iterations and the horizon length sufficiency tests, though we may ultimately solve for problem horizons up to $2N_{M}$, as indicated by \eqref{ucpadp-approx-horzlen}. Ultimately this choice is arbitrary, and it is possible for some problems to use another scheme for selecting a new length.

From a practical perspective, we have found that setting $\varepsilon_{\mu}\approx2d_{u}$ and $\varepsilon_{x}\approx2d_{x}$ (the distance between points in $\mathcal{U}$ and $\mathcal{X}$ respectively) is a good design choice for well-behaved problems. Smaller values raise the risk of never terminating, e.g.~due to residual state trajectory jitter caused by approximation inherent to interpolation, while larger values give an unnecessarily large approximation of the true control policy $\mu_{R,0}^{*}$. Ultimately, this implies that UCPADP has to some degree only one tuning parameter: the ADP discretisation, which trades off accuracy with computational time and memory demands.

As UCPADP is based on interpolating ADP (and in turn DP) it is subject to the inherent limitations of DP methods, in particular its poor scaling with problem dimensionality (colloquially referred to as the ``curse of dimensionality'') \citep{bellman1954theory,Bertsekas2017}. This limits UCPADP to low- to moderate-dimensional problems. The examples shown in the following section (with two state variables and one control variable, giving a total of three independent variables) are easily solved using an ordinary desktop computer on the order of one minute to one hour (depending on the demanded solution accuracy). In practice, we expect UCPADP to be viable for up to 4--6 continuous-variable problems, depending on the discretisation of the state and control variables, the nature of the problem, and the available computational power.

A general implementation of the UCPADP method in the MATLAB language, including the numerical examples in the following section, is available at \texttt{\url{https://gitlab.com/lerneaen_hydra/ucpadp}}.

\section{Representative examples\label{sec:Representative-examples}}

We illustrate the UCPADP method, introduced in \secref{The-ucpadp-method}, by solving two simple problems. Though ``toy'' problems in some sense, recall that \assmref{base-assumptions} allows for significantly more difficult (and practically relevant) problems. First we consider the classical minimum-time inverted pendulum problem, where we highlight the stopping criterion of UCPADP. Afterwards, we consider the problem of maintaining an average pendulum angle with minimum control power, illustrating the average-constraint properties shown in \subsecref{Average-constraint-relaxation}.

\begin{figure}
\begin{centering}
\includegraphics[width=0.4\columnwidth]{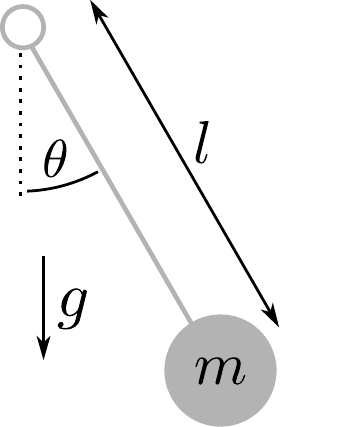}
\par\end{centering}
\caption{A simple pendulum. \label{fig:A-simple-pendulum}}

\end{figure}

We will consider the dynamical system given by a simple pendulum (\figref{A-simple-pendulum}) for both problems. For a pendulum with length $l$, point mass $m$, gravitational force $g$, damping coefficient $d$, angle $\theta$, and applied torque $u$, the dynamic equation for the system can be derived as
\begin{equation}
\ddot{\theta}+\frac{d}{m}\dot{\theta}+\frac{g}{l}\sin\left(\theta\right)=\frac{1}{ml^{2}}u.\label{eq:pendulum-dyn}
\end{equation}
In both the following examples we will assume a discrete-time control system with sample rate $t_{s}$, i.e.\ the control input $u$ is piecewise constant over intervals of uniform time $t_{s}$. If the problem is reformulated as a set of coupled first-order ordinary differential equations with a state variable vector
\begin{align}
x & \triangleq\left[\begin{array}{c}
\theta\\
\dot{\theta}
\end{array}\right]
\end{align}
then we can express the state at the next sample as
\begin{equation}
x_{k+1}=f_{p}\left(x_{k},u_{k}\right),
\end{equation}
where $f_{p}$ is given by solving \eqref{pendulum-dyn} over a time $t_{s}$ with initial condition $x_{k}$ and constant control input $u_{k}$.

\subsection{The inverted pendulum\label{subsec:The-inverted-pendulum}}

To illustrate the mechanics of UCPADP's termination criterion, consider the traditional minimum-time inverted pendulum problem (formulated here as an infinite-horizon problem)
\begin{subequations}
\label{eq:min-time-inverted}
\begin{align}
J^{*} & =\min_{\bar{x},\bar{u}}\lim_{N\to\infty}\frac{1}{N}\sum_{k=0}^{N-1}f_{c}\left(x_{k}\right)\\
f_{c}\left(x\right) & =\begin{cases}
0 & \text{if }|\theta_{k}-\pi|<2d_{x},|\dot{\theta}_{k}|<2d_{x}\\
1 & \text{else}
\end{cases}\label{eq:min-time-inverted-cost}\\
\shortintertext{\text{subject to}}x_{k+1}= & f_{p}\left(x_{k},u_{k}\right)\\
\left|u_{k}\right|\leq & 1\label{eq:min-time-inverted-control}\\
-2\leq\theta_{k} & \leq3.5\label{eq:min-time-inverted-theta}\\
-1.5\leq\dot{\theta}_{k} & \leq2.\label{eq:min-time-inverted-thetadot}
\end{align}
\end{subequations}

All the following results are shown for a sample time of $t_{s}=0.2$, pendulum parameters set to give the system dynamics equation $\ddot{\theta}+\sin\left(\theta\right)=u$, state variables discretised by a Cartesian grid with separation $d_{x}=[0.05,0.05]^{T}$ in the range allowed by \eqref{min-time-inverted-theta} and \eqref{min-time-inverted-thetadot}, and the control variable discretised with even spacing $d_{u}=0.01$ in the range allowed by \eqref{min-time-inverted-control}. Setting $\varepsilon_{x}$ and $\varepsilon_{\mu}$ to the suggested value of twice the discretisation gives $\varepsilon_{x}=\left[0.1,0.1\right]^{T}$ and $\varepsilon_{\mu}=0.02$. 

Note that the cost function \eqref{min-time-inverted-cost} equally penalises all pendulum configurations other than the single vertical zero-velocity state combination, and with an infinite horizon (and small enough $d_{x}$) gives a solution arbitrarily close to the traditional minimum-time formulation. The state bounds \eqref{min-time-inverted-theta} and \eqref{min-time-inverted-thetadot} have been chosen to give a reasonable range for the specific initial value we will study shortly.

\begin{figure}
\begin{centering}
\includegraphics[width=0.66\columnwidth]{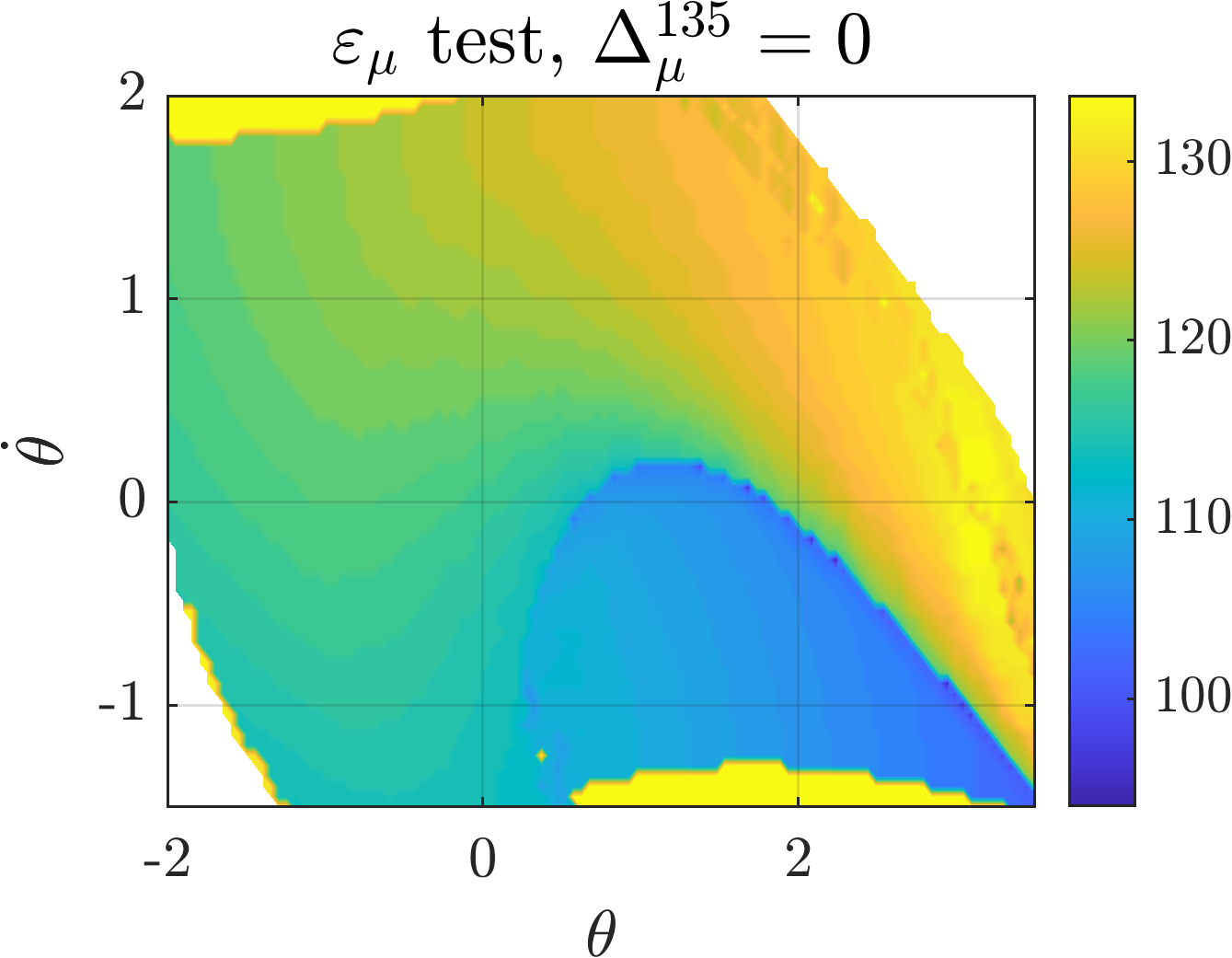}
\par\end{centering}
\caption{Visualisation of $\varepsilon_{\mu}=0.02$ test for $N_{M}'=135$. The coloured regions indicate the number of samples that the control policy varies less than $\varepsilon_{\mu}$, while white regions indicate a feasible solution could not be found, i.e.~white regions lie outside of $\mathcal{F}'_{135}$. Note that as $\mathcal{U}$ is discrete then $\Delta_{\mu}^{135}$ is also, i.e.~we here have $\Delta_{\mu}^{135}$ identically equal to zero. \label{fig:ucpadp_termination_eps_u}}
\end{figure}

\begin{figure}
\begin{centering}
\includegraphics[width=0.66\columnwidth]{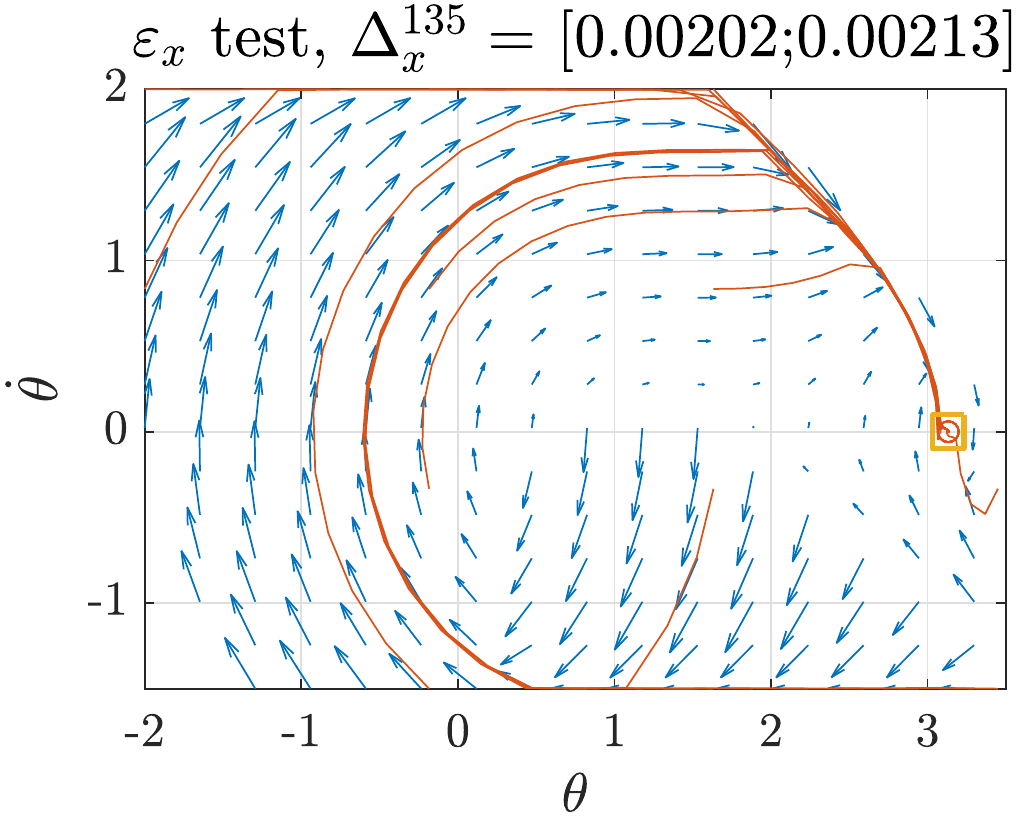}
\par\end{centering}
\caption{Visualisation of $\varepsilon_{x}=\left[0.1,0.1\right]^{T}$ test for $N_{M}'=135$. Blue arrows indicate the motion of the system through its phase space. Representative closed-loop trajectories are shown in red. The closed-loop state at $N=135$ is shown by (overlapping) small red circles near $\theta=\pi$, $\dot{\theta}=0$. The yellow box indicates the $\Delta_{x}^{135}$ termination criterion, which is satisfied as all states at $N=135$ lie inside the box. \label{fig:ucpadp_termination_eps_x}}
\end{figure}

\begin{figure}
\begin{centering}
\includegraphics[width=0.66\columnwidth]{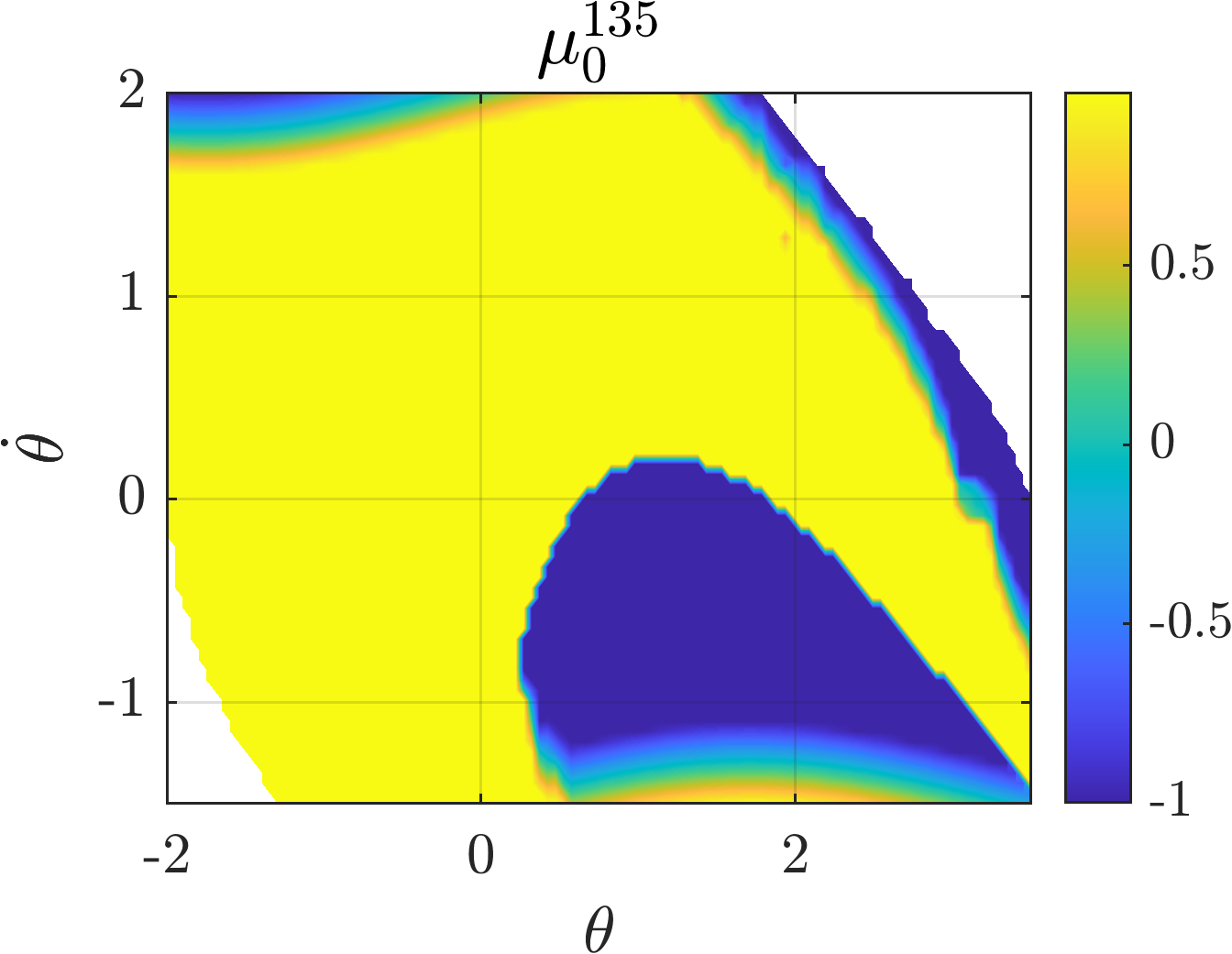}
\par\end{centering}
\caption{Control policy associated with \eqref{min-time-inverted} for horizon $N=135$. The coloured region shows the optimal control to apply for any given state, while white regions indicate infeasible states, i.e.~outside of $\mathcal{F}'_{135}$. \label{fig:inverted-pendulum-control-law}}
\end{figure}

For the above problem, UCPADP terminates after testing a horizon of $N_{M}'=135$, indicating that $45<N_{M}\leq135$. An illustration of termination criterion \eqref{cond_approx_min_stationary_eps_mu} is shown in \figref{ucpadp_termination_eps_u}, where we can verify the condition is satisfied as all values are above $\left\lceil N/2\right\rceil =68$. Furthermore, $\Delta_{\mu}^{k}$ will by construction take values from $\mathcal{U}=\left\{ 0,\pm0.01,\pm0.02,\dots,\pm1\right\} $. For $\varepsilon_{\mu}=0.02$ \eqref{cond_approx_min_stationary_eps_mu} will thus only be satisfied for values $-0.01,0,0.01$. We can see this in \figref{ucpadp_termination_eps_u}, where $\Delta_{\mu}^{135}=0$. Similarly, criterion \eqref{cond_approx_min_stationary_eps_x} is illustrated in \figref{ucpadp_termination_eps_x}, where we can verify that representative trajectories all converge to a region bounded by $\varepsilon_{x}$ (shown by the yellow box). An illustration of the control policy ultimately generated by UCPADP is shown in \figref{inverted-pendulum-control-law}. Solving this specific problem took approximately 10 minutes using a standard desktop PC.

\begin{figure}
\begin{centering}
\includegraphics[width=1\columnwidth]{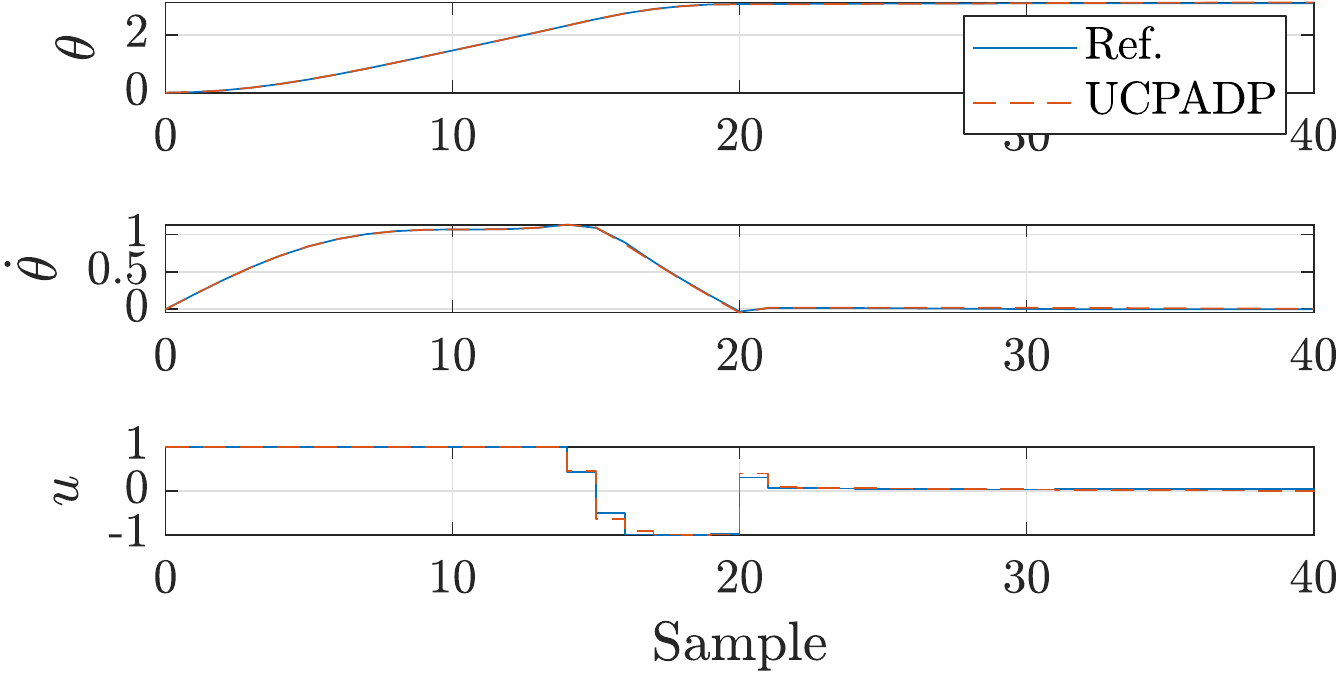}
\par\end{centering}
{\footnotesize{}\caption{Comparison of the solutions given by UCPADP and an open-loop ADP reference method for $x_{0}=\left[0,0\right]^{T}$.\label{fig:ucpadp-vs-ref-invert}}
}{\footnotesize\par}
\end{figure}

\figref{ucpadp-vs-ref-invert} shows a comparison of the solution generated by UCPADP and a reference solution, generated by formulating a problem with an explicit horizon of $N=10N_{M}'=1350$ (i.e.~one order of magnitude longer the UCPADP horizon), for $x_{0}=\left[0,0\right]^{T}$. Here the reference solution is generated using a traditional ADP scheme, configured with the same sample time and state/control grid discretisation. Note that we intentionally compare the UCPADP solution to a traditional ADP solution (in contrast to, for instance, an analytical solution) as we wish to highlight the accuracy of the automatically sized horizon, rather than the accuracy of an interpolating ADP scheme.

The average cost over the time interval shown in \figref{ucpadp-vs-ref-invert} is $0.09664$ for the UCPADP solution, while the cost associated with the reference solution is $0.09689$, i.e.~a deviation\footnote{The fact that the UCPADP solution has a lower associated cost is likely due to the inherent approximation of interpolating ADP.} of 0.25\%. We can conclude (for this specific problem) that the cost associated with the UCPADP solution is virtually identical to a conventional ADP solution, indicating that the identified horizon $N_{M}'=135$ was sufficient. 

\subsection{The constant-angle pendulum}

Let us now consider a problem that illustrates the properties of the average constraint introduced in \subsecref{Average-constraint-relaxation}. Assume we wish to solve
\begin{subequations}
\label{eq:angled-pendulum}
\begin{align}
J^{*} & =\min_{\bar{x},\bar{u}}\lim_{N\to\infty}\frac{1}{N}\sum_{k=0}^{N-1}u_{k}^{2}\label{eq:angled_pendulum_unaug_cost}\\
\shortintertext{\text{subject to}}x_{k+1} & =f_{p}\left(x_{k},u_{k}\right)\\
\lim_{N\to\infty}\frac{1}{N}\sum_{k=0}^{N-1}\theta_{k} & =\theta_{\text{ref}}\label{eq:angled_pendulum_avg_const}\\
|u_{k}| & \leq1,|\theta_{k}|\leq1,|\dot{\theta}_{k}|\leq1,
\end{align}
\end{subequations}
i.e.\ the problem of keeping the average pendulum angle at a setpoint $\theta_{\text{ref}}$ while minimising the quadratic control input $u_{k}^{2}$.

By \thmref{relaxed} we can avoid including the average constraint \eqref{angled_pendulum_avg_const} by augmenting the cost functional \eqref{angled_pendulum_unaug_cost} as
\begin{equation}
J_{R}^{*}=\min_{\bar{x}_{R},\bar{u}_{R}}\lim_{N\to\infty}\frac{1}{N}\sum_{k=0}^{N-1}u_{k}^{2}+\lambda\theta_{k}\label{eq:angled_pendulum_aug_cost}
\end{equation}
for a constant scalar $\lambda$. Assuming the problem reaches an equilibrium with control $u_{\text{eq}}$ and states $\theta=\theta_{\text{ref}},\dot{\theta}=0$, by \eqref{pendulum-dyn} we have $u_{\text{eq}}=mgl\sin(\theta_{\text{ref}})$. We can thus express the equilibrium cost as 
\begin{equation}
c_{\text{eq}}=(mgl\sin(\theta_{\text{ref}}))^{2}+\lambda\theta_{\text{ref}}\label{eq:eq_cost_reform}
\end{equation}

which is a function of one variable. Equation~\eqref{eq_cost_reform} has one unique stationary point (a minimum) in the permissible range $|\theta|<1$, and we can thus find the specific value $\lambda$ that gives the lowest equilibrium cost at the desired setpoint by setting $\frac{dc_{\text{eq}}}{d\theta_{\text{ref}}}=0$ and solving for $\lambda$, giving
\begin{equation}
\lambda_{0}=-2m^{2}g^{2}l^{2}\sin\left(\theta_{\text{ref}}\right)\cos\left(\theta_{\text{ref}}\right).
\end{equation}

We can now reformulate \eqref{angled-pendulum} as the equivalent problem
\begin{subequations}
\label{eq:angled-pendulum-reform}
\begin{align}
J_{R}^{*} & =\min_{\bar{x}_{R},\bar{u}_{R}}\lim_{N\to\infty}\frac{1}{N}\sum_{k=0}^{N-1}u_{k}^{2}-\theta_{k}\lambda_{0}\label{eq:angled-aug-cost}\\
\shortintertext{\text{subject to}} & x_{k+1}=f_{p}\left(x_{k},u_{k}\right),\left|u_{k}\right|\leq1,|\theta_{k}|\leq1,|\dot{\theta}_{k}|\leq1.
\end{align}
\end{subequations}

\begin{figure}[h]
\begin{centering}
\includegraphics[width=1\columnwidth]{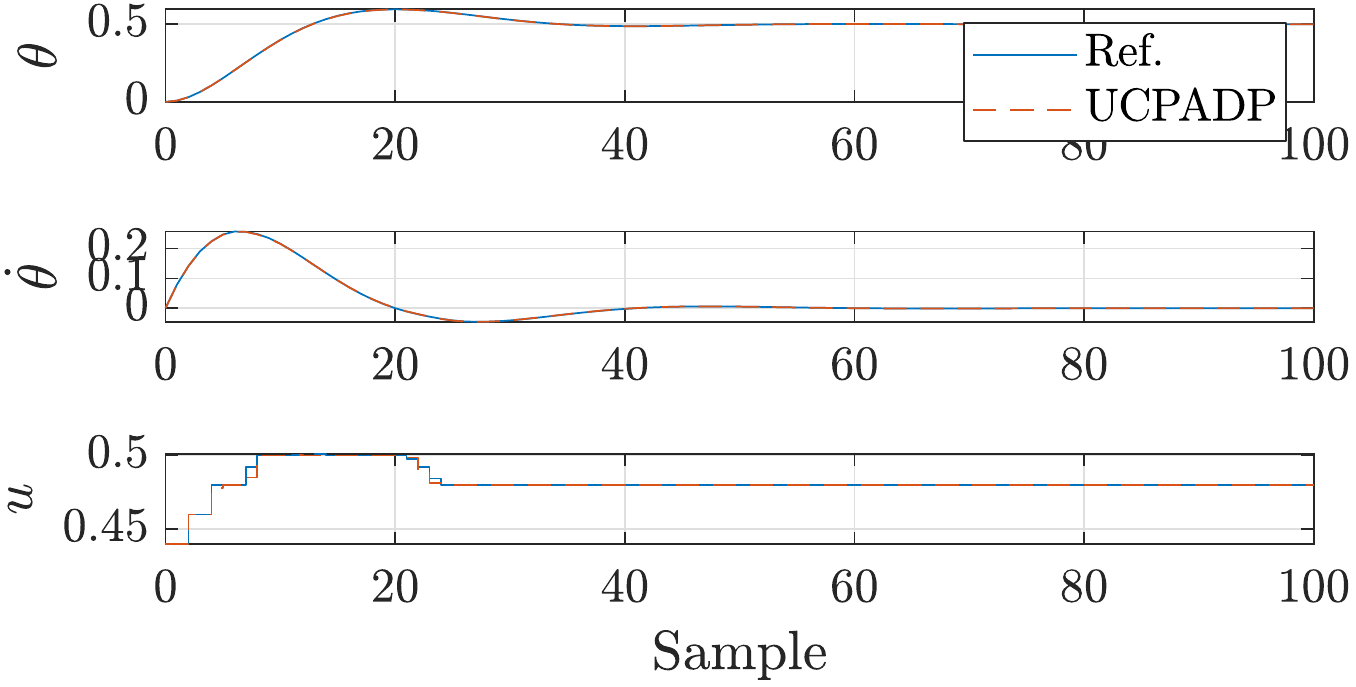}
\par\end{centering}
\caption{Comparison of the solutions given by UCPADP and an open-loop ADP reference method for $x_{0}=\left[0,0\right]^{T}$. The state trajectories are nearly identical and the control trajectory displays only very small differences. \label{fig:ucpadp-vs-ref-angled}}
\end{figure}

As in the previous example, we discretise the state and control variables evenly in the permissible space, here with separation $d_{x}=[0.02,0.02]^{T}$ and $d_{u}=0.01$ respectively. Solving \eqref{angled-pendulum-reform} for pendulum parameters resulting in a system dynamics equation $\ddot{\theta}+\dot{\theta}+\sin\left(\theta\right)=u$ and $\theta_{\text{ref}}=0.5$ gives the results shown in \figref{ucpadp-vs-ref-angled} (again compared with a reference solution given by explicitly choosing a large horizon, one order of magnitude larger than the horizon given by UCPADP).

\begin{figure}
\begin{centering}
\includegraphics[width=1\columnwidth]{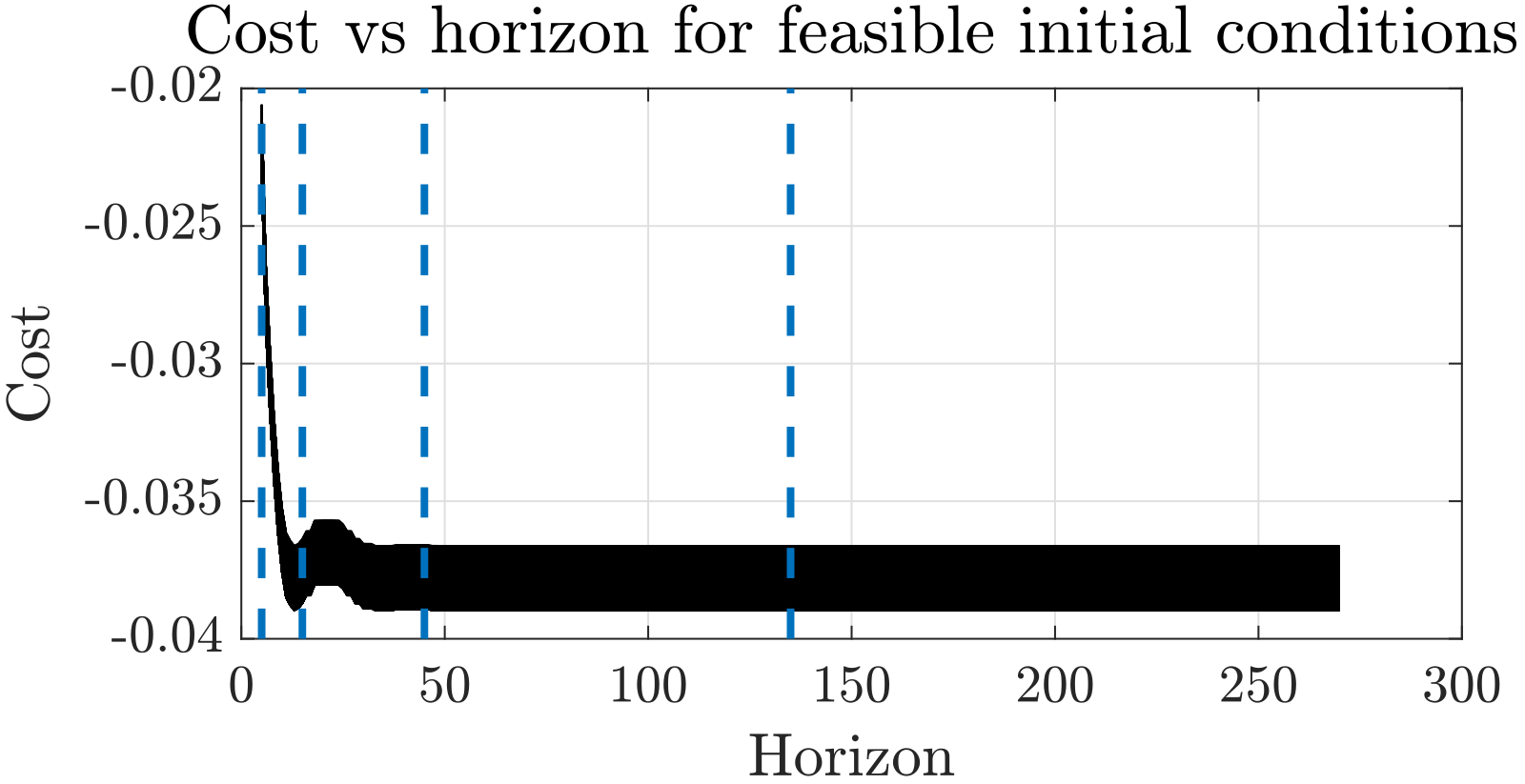}
\par\end{centering}
\caption{Cost of the applying control policy associated with varying problem horizon lengths. The cost is shown for all feasible initial conditions, resulting in a range of costs for any given horizon (e.g.~$\left[-0.037,-0.039\right]$ for horizons $\protect\geq40$). Problem horizons tested by UCPADP are shown with dashed lines.\label{fig:cost-vs-horizon-ensemble}}
\end{figure}

For this problem, we find that the UCPADP solution generates a solution with control cost (i.e.~$\sum u_{k}^{2}$) of $0.2321$ over the horizon shown in \figref{ucpadp-vs-ref-angled}, while the control cost associated with the reference solution is $0.2323$ (i.e.~a deviation of 0.09\%), again showing that the accuracy of the UCPADP solution is virtually identical to that of a reference ADP solution. 

For comparison, in \figref{cost-vs-horizon-ensemble} we also show the solution quality parameterised by different finite horizons. More specifically, we solve the finite-horizon counterpart of \eqref{angled-pendulum-reform}, i.e.~using the notation introduced in \eqref{finite-horz-prob-def}, for varying finite horizons $N$ (denoted the \emph{problem horizon}), resulting in the associated control policies $\mu_{R,0}^{*N}$. We then apply the control policy to the set of initial conditions feasible with a long horizon $N=1350$ (denoted the \emph{trajectory horizon}). The plot shows the augmented cost of the trajectory horizons, i.e.~$J_{R}^{*}$, parameterised by different problem horizons. We can identify that the average cost is higher for short problem horizons than for long problem horizons, and that the cost associated with problem horizons $\gtrsim40$ is constant, indicating that for this problem a problem horizon $\gtrsim40$ is sufficient. It may therefore seem like UCPADP is inefficient in its choice of problem horizon (135 samples). However, computing the average cost of any given problem horizon shown in \figref{cost-vs-horizon-ensemble} is time consuming, with each individual problem horizon taking approximately the same time to compute as the \emph{entire} UCPADP solution, as well as requiring problem-specific knowledge of the initial conditions and trajectory horizon to average over. The trade-off between spending time computing additional back-calculation steps and checking whether a given horizon is sufficiently large thus motivates a scheme like our proposed horizon-doubling method.

\section{Conclusions}

In this paper we have introduced UCPADP, a numerical method inspired by API. UCPADP can be used to generate a near-optimal control policy for general undiscounted continuous-valued infinite-horizon nonlinear optimal control problems. The problem can also optionally be constrained to converge to a given equilibrium. The primary contribution of UCPADP is the introduction of a termination criterion that is amenable to the undiscounted case, while still allowing for general costs and constraints. We have evaluated the method by solving two simple, but representative, problems. For both examples we showed that the generated control policy was on par with the accuracy of a reference ADP solution (whose accuracy is determined by the chosen discretisation of the problem).

UCPADP has several properties that render it useful as as one part of the process of constructing an on-line controller. Firstly, it shares a property with other API methods in that it does not require any a-priori information about a suitable horizon, instead performing an indefinite number of iterations and terminating when a suitable problem horizon is found. Secondly, the tuning parameters are simple to grasp, as they trade off solution accuracy with computational time and memory demands. Finally, the output from UCPADP, as with other API methods, is a control policy (i.e.~a state feedback table). After this control policy is computed in an off-line phase it can in turn be used to construct a subsequent on-line controller with very low computational demand, only requiring a simple interpolation operation to determine the control signal.

Full source code of the implementation as well as the specific problems studied is available at \texttt{\url{https://gitlab.com/lerneaen_hydra/ucpadp}.}

\section*{Funding}
This work was performed within the Combustion Engine Research Center at Chalmers (CERC) with financial support from the Swedish Energy Agency.

\bibliographystyle{apacite}
\bibliography{ref}

\end{document}